# Identification of Feasible Regions Using R-Functions

S. Kucherenko[a*], N. Shah[a], O.V. Klymenko[b*]

[a]Imperial College London, London, SWT 2AZ, UK

s.kucherenko@imperial.ac.uk

[b]School of Chemistry and Chemical Engineering, University of Surrey, Guildford GU2 7XH, UK

Abstract.

The primary objective of flexibility analysis is to identify and define the feasibility region, which represents the range of operational conditions (e.g., variations in process parameters) that ensure safe, reliable, and feasible process performance. This work introduces a novel flexibility analysis method that requires only that model constraints (e.g., defining product Critical Quality Attributes or process Key Performance Indicators) be explicitly provided or approximated by a closed-form function, such as a multivariate polynomial model. The method is based on V.L. Rvachev's R-functions, enabling an explicit analytical representation of the feasibility region without relying on complex optimization-based approaches. R-functions offer a framework for describing intricate geometric shapes and performing operations on them using implicit functions and inequality constraints. The theory of R-functions facilitates the identification of feasibility regions through algebraic manipulation, making it a more practical alternative to traditional optimization-based methods. The effectiveness of the proposed approach is demonstrated using a suite of well-known test cases from the literature.

1. Introduction

Flexibility analysis, introduced by Halemane and Grossmann (1983) and Grossmann and Morari (1983), plays a crucial role in process design and operation by determining the conditions under which a system remains feasible under specified constraints, such as Critical Quality Attributes (CQAs).

Two classical flexibility analysis problems are the flexibility test and the flexibility index (Swaney and Grossmann, 1985a, 1985b). The flexibility test determines whether feasible operation can be achieved within a given range of uncertainty scenarios. The flexibility index, on the other hand, quantifies the operational range, representing the maximum scaled deviation of all process parameters from nominal conditions. This can be visualized as the largest rectangle inscribed within the feasible space, where steady-state operation can be maintained by adjusting control variables.

The methodology developed by Grossmann and colleagues has been widely recognized as a powerful tool for identifying uncertainty ranges within which a process or product can operate feasibly. It has found applications in various fields, including chemical engineering, pharmaceutical manufacturing, and energy systems, where managing uncertainty is critical for ensuring process reliability and product quality.

Dimitriadis and Pistikopoulos (1995) extended flexibility analysis to dynamic systems, incorporating time-dependent uncertainties. Their work distinguished dynamic flexibility analysis from the earlier steady-state approaches, enabling the evaluation of processes subject to transient or time-varying conditions.

Ierapetritou (2001) noted that existing methods tended to significantly underestimate the actual feasible region. Consequently, developing a more accurate determination of the feasible range is of great importance for analysing uncertainties in process and product design. Goyal and Ierapetritou (2002) introduced a framework for accurate approximation of the feasible region or operating envelope of a given design based on the basic idea of approximating the feasible region from the inside using the simplicial approximation approach and from the outside using the tangent planes at specific boundary points.

Banerjee and Ierapetritou (2005) noticed that feasibility analysis shares similarities with surface reconstruction, as both focus on identifying and accurately estimating the boundary of a given region. Traditional approaches to feasibility analysis approximate this boundary using linear inequalities, either through the use of a hyperrectangle or by approximating the convex hull within the feasible space. These methods perform well only for convex, connected regions, but they are inaccurate when dealing with nonconvex or disjoint regions. In contrast, surface reconstruction techniques can effectively capture both nonconvex and disjoint regions by defining the bounding surface through piecewise linear functions. To approximate the region using the $\alpha$-shape technique, the first step is to generate points representing the region by sampling the feasible space. Traditional sampling methods assume a uniform distribution, resulting in uniform coverage of the parameter space without distinguishing between feasible and infeasible regions. To address this, Banerjee and Ierapetritou (2005) proposed a new sampling method that focuses on the smaller feasible region by reformulating the sampling process as an optimization problem solved using a genetic algorithm (GA).

Once the feasible region is approximated as a polygon, point-in-polygon tests determine if a point lies inside. Banerjee and Ierapetritou (2005) proposed to use the ray tracing algorithm. With $\mathcal{O}(n)$, complexity, this method becomes inefficient for large polygons or multiple points and struggles with edge cases. While it supports non-convex shapes, it is not easily extendable to 3D polyhedra, requiring more complex computations.

Boukouvala and Ierapetritou (2012) further advanced feasibility evaluation by using Kriging interpolation and developing an adaptive sampling strategy to reduce sampling costs while preserving the accuracy of the feasibility space. Grossmann et al. (2014) review the early mathematical formulations and solution methods developed by Grossmann and colleagues for quantifying Static Resiliency (Flexibility), along with a discussion of recent advancements in the field since these foundational studies.

Although optimization methods used in flexibility analysis are computationally efficient, they may lead to non-unique solutions and complex multi-level problems. Additionally, these methods are typically limited to processes with closed-form and differentiable constraints.

Determining the feasible region in problems with black-box constraints can be challenging, when these constraints arise from empirical models, or closed-source software. Black-box constraint functions may not be explicitly differentiable, hence efficient gradient-based optimization methods may not be practical for feasibility analysis. In such cases surrogate-based approaches can be used to approximate the feasibility function and identify the feasible operating region (Banerjee et al., 2010; Boukouvala and Ierapetritou, 2012).

Zhao and Chen (2018) proposed an operational flexibility analysis method to explicitly describe the feasibility region and establish functional relationships between uncertain parameters and control variables. The problem is framed as an existential quantifier problem, solved using the cylindrical algebraic decomposition method to eliminate quantifiers and derive triangular quantifier-free formulas. A logical combination of these formulas defines the feasibility region and provides a control strategy for feasible operations. Although it is limited to polynomial functions, the method can handle non-polynomial systems using local polynomial approximations. Due to the computational intensity of symbolic computation needed to apply the cylindrical algebraic decomposition method, the computational burden increases exponentially with system size, hence the method is only applicable to relatively small-scale problems.

One key application of feasibility analysis is in the field of pharmaceutical processes. The pharmaceutical manufacturing industry is subject to strict requirements to ensure drug quality and safety. This drives the Quality-by-Design (QbD) approach adopted by regulatory agencies worldwide. QbD focuses on robust process design to maintain product quality, with its core concept being the design space (DS) defined as the multivariate region where acceptable-quality material can be produced (ICH, 2009). In pharmaceutical manufacturing, feasibility analysis helps define this space by identifying the feasible operating region where all quality constraints are met (Boukouvala et al., 2010a).

Identification of probabilistic DS (Castagnoli et al., 2010) is a demanding task and for a typical practical problem the traditional approach based on exhaustive sampling requires costly computations. A novel theoretical and numerical framework for determining probabilistic DS using metamodelling and adaptive sampling was proposed in Kucherenko et al. (2020). It is based on the multi-step adaptive technique using a metamodel for a probability map as an acceptance-rejection criterion to optimize sampling for identification of the DS.

Flexibility index and design centering are two complementary approaches for evaluating candidate DSs. From a geometric viewpoint, the flexibility index seeks to identify the largest hyperrectangle (or other geometric shape) that can be inscribed within the feasible region, while design centering aims to locate the center of that region, maximizing robustness to variations. Zhao et al. (2022) proposed a novel flexibility index formulation using a direction search method, capable of handling design spaces of any shape. However, its application was demonstrated only for rectangular and elliptical cases, highlighting the method's complexity. Tian et al. (2024) presented an overview of key concepts and computational methods for evaluating feasibility and flexibility, including recent advances and applications in optimizing process design and operations.

Flexibility analysis and DS identification involve defining a region—whether a feasible space or a design space—where specific constraints are satisfied while accounting for uncertainties and ensuring robustness against input variability. The methodologies used in DS analysis, can be adapted to characterize feasibility regions in pharmaceutical or chemical processes.

Kucherenko et al (2024) introduced the R-DS method, which utilizes V.L. Rvachev's R-function methodology (Rvachev, 1982) to anaalytically represent the DS. By leveraging R-functions, complex geometric shapes can be represented and manipulated through implicit functions, enabling the identification of the DS via algebraic operations. Unlike previous studies where the DS was identified by simultaneously considering all CQAs, the developed approach involves a two-stage process. First, feasible subspaces of the process operational parameter space are identified for each constraint separately, deriving implicit functions to analytically describe their boundaries. Next, the theory of R-functions is used to define the DS as the intersection of these feasible subspaces. This method yields the DS given by a single analytical expression. Moreover, the analytical definition of the DS using R-functions offers additional advantages, such as streamlining operations like "boxing" (determining a tight axis-aligned bounding box), identifying points within the feasible domain, and calculating the domain's volume. In addition, this method also allows the handling of non-convex and disjoint DS shapes.

In this study, we demonstrate the effectiveness of V.L. Rvachev's R-functions in identifying and characterizing feasible regions for process design and operational flexibility analysis. While all the test cases considered involve constraints with closed-form equations, when constraints are not defined explicitly, the methodology from our previous work on the analytical identification of DS can be applied (Kucherenko et al., 2020). This approach involves approximating constraints with closed-form functions, such as multivariate polynomial models.

This paper is organized as follows. Section 2 presents the mathematical formulation of the feasibility analysis problem, defining the feasibility function. Section 3 provides an overview of R-functions, including their theoretical background and application in describing geometric regions. Section 4 introduces the methodology for identifying the feasible region using R-functions. Section 5 presents test cases from the literature, demonstrating the application of R-functions to feasibility analysis in both open-loop and closed-loop scenarios. Section 6 discusses the advantages of the R-function methodology, comparing it with traditional methods such as simplicial approximation and α-shape reconstruction. Section 7 provides conclusions, highlighting the benefits of the proposed approach and discussing potential extensions for handling complex constraint functions.

## 2. Feasibility analysis. Mathematical formulation

Feasibility analysis can be mathematically formulated as the maximum of the process constraint violation:

$$\psi(d, x) = \min_{z} \max_{j \in J} \{g_j(d, z, x)\} \text{ s.t. } x \in T, \quad (1)$$

where $\psi(d, x)$ is the process feasibility function; $d$ is the vector of the process design variables, $z$ is the vector of process variables that can be manipulated, and $x$ is the vector of critical input variables (e.g., critical process parameters and material properties); $g_j$ are the functions of the $J$ process constraints in the form $g_j(d, z, x) \leq 0$ which must be satisfied during process operation; $T$ describes the range in which uncertain inputs can vary:

$$T = \{x \mid x^L \leq x \leq x^U\},$$

where $x^L$ and $x^U$ are lower and upper bounds, respectively.

Assuming that design variables $d$ are constant, Eq. (1) simplifies as follows:

$$\psi(x) = \min_{z} \max_{j \in J} \{g_j(z, x)\} \text{ s.t. } x \in T. \quad (2)$$

The value of the feasibility function in Eq. (2) indicates whether the process is feasible. The function $\psi(x)$ measures the "distance" to feasibility.

- If $\psi(x) > 0$, there is at least one constraint that is being violated. The process is not feasible.
- If $\psi(x) \leq 0$ no constraint is being violated. The process is feasible. When $\psi(x) = 0$, vector $x$ identifies the boundary of the feasible region.

When $\psi(x) = 0$, the system is exactly feasible.

In order to identify the process DS, we need to determine the values of each variable $x_i$ in vector $x$, while adjusting the values of each manipulated variable $z_k$ in vector $z$, so that the condition $\psi(x) \leq 0$ is satisfied. Values of control actions $z_k$ in vector $z$ minimizing the maximum of process constraints, $\max_{j \in J}\{g_j(z, x)\}$, have been traditionally computed through nonlinear sequential quadratic programming (SQP) optimization.

Problem (2) can be posed as a standard optimization problem (LP or NLP) by defining a scalar variable $u$, such that:

$$\psi(d, x) = \min_{z, u} u$$
$$\text{s.t. } g_j(d, z, x) \leq u, \forall j \in J. \quad (3)$$

## 3. Theory of R-Functions

A function $y = R(x_1, \ldots, x_n)$ is called an R-function it there exists a Boolean function $Y = F(X_1, \ldots, X_n)$ such that the following relationship holds for any $(x_1, \ldots, x_n) \in \mathbb{R}^n$ (Shapiro, 1991):

$$S[R(x_1, \ldots, x_n)] - F[S(x_1), \ldots, S(x_n)]. \tag{4}$$

Here, $S(x)$ is a step function defined as:

$$S(x) = \begin{cases} 0, & \text{if } x \leq 0 \\ 1, & \text{if } x \geq 0 \end{cases}$$

The Boolean variables $X_1, \ldots, X_n$ are derived by applying the step function $S(x)$ to each real variable: $X_i = S(x_i)$ for $i = 1, 2, \ldots, n$. A Boolean function is a mathematical function that maps binary inputs (typically 0 and 1, or false and true) to a binary output. The Boolean function $F$ is referred to as the companion function of the R-function.

A real-valued function is classified as an R-function if its sign changes only when one or more of its arguments change sign. Basic R-functions can be constructed using logical operations on implicit functions:

- Union (OR): The union of two implicit functions $R_1$ and $R_2$ is defined by $R_\vee \geq 0$ if $R_1 \geq 0$ or $R_2 \geq 0$. This is also known as R-disjunction.
- Intersection (AND): The intersection of two implicit functions $R_1$ and $R_2$ is defined by $R_\lambda \geq 0$ if $R_1 \geq 0$ and $R_\wedge \geq 0$. This is also called R-conjunction.
- Complement (NOT): The negation of $R$ describes the complement of the set defined by this function: $R_\neg \geq 0$ if $R \leq 0$. This is also referred to as R-negation.

In all the above definitions, the functions depend on the input $(x_1, \ldots, x_n)$. These operations use specific algebraic combinations of implicit functions to preserve geometric properties. According to the above definition, the following real-valued functions qualify as R-functions, with their respective Boolean companion functions indicated in parentheses (Shapiro, 1991):

$$\begin{aligned} C &\equiv \text{const}; & \text{(logical 1)} \\ \bar{x} &\equiv -x; & \text{(logical negation } \neg \text{)} \\ x_1 \wedge_1 x_2 &\equiv \min(x_1, x_2); & \text{(logical conjunction } \wedge \text{)} \\ x_1 \vee_1 x_2 &\equiv \max(x_1, x_2); & \text{(logical disjunction } \vee \text{)} \end{aligned} \tag{5}$$

While R-functions are infinite in variety, only a sufficiently complete system is needed for practical applications. For example, conjunction ($\wedge$) and negation ($\neg$) form a complete system. In this work, we consider the following general set $R_\alpha$ of R-functions:

$$\begin{aligned} x_1 \wedge_\alpha x_2 &\equiv \frac{1}{1+\alpha}\left(x_1 + x_2 - \sqrt{x_1^2 + x_2^2 - 2\alpha x_1 x_2}\right), \\ x_1 \vee_\alpha x_2 &\equiv \frac{1}{1+\alpha}\left(x_1 + x_2 + \sqrt{x_1^2 + x_2^2 - 2\alpha x_1 x_2}\right), \end{aligned} \tag{6}$$

where $\alpha(x_1, x_2)$ is an arbitrary symmetric function such that $-1 < \alpha(x_1, x_2) \leq 1$. The value of $\alpha$ determines the smoothness of the R-functions. The first function in (6) is called R-conjunction, while the second is R-disjunction. Negation $\neg$ is defined as

$$\bar{x} \equiv -x. \tag{7}$$

In Euclidean space $\mathbb{R}^n$, the set of points where the equality $R(x_1, \ldots, x_n) = 0$ holds is called a drawing, while the set of points satisfying $R(x_1, \ldots, x_n) \geq 0$ is referred to as a region. Many familiar geometric objects, such as curves and surfaces, are algebraic drawings, often defined by polynomial functions. These functions act as fundamental building blocks for constructing more complex geometric shapes of practical interest. Such constructions are achieved through standard set operations or, equivalently, by combining algebraic primitives using logical operations.

A composite geometric region $D \subseteq \mathbb{R}^n$ can be described as
$$D = F[(\varphi_1 \geq 0), \ldots, (\varphi_m \geq 0)], \tag{8}$$
where real-valued inequalities $\varphi_i(x_1, \ldots, x_n) \geq 0$ define primitive geometric regions, and $F$ is a set function constructed using standard set operations $\cap, \cup, \setminus$ on these primitive regions. The function $F$ can be viewed as a Boolean function, in which the set operations are replaced with the corresponding logical functions $\wedge, \vee, \neg$.

In practice, we seek a single real-function inequality $R(x_1, \ldots, x_n) \geq 0$ that defines the composite object $D$ as a closed subset of $\mathbb{R}^n$. It has been proven (Rvachev, 1982; Shapiro, 1991) that to obtain a real function defining the region $D$, it suffices to construct an appropriate R-function and substitute its arguments with the real functions $\varphi_i$ defining the primitive regions:
$$R(\varphi_1, \ldots, \varphi_m) \geq 0. \tag{9}$$

4. Application of R-functions in feasibility analysis

In this section, we use the theory of R-functions to reformulate the feasibility test problem defined in Eq. (2) or (3) above in terms of finding an explicit analytical representation of the feasibility function. We name this approach R-feasibility for brevity.

4.1. Open-loop scenario

In the case of an open-loop scenario, in which the system operates without feedback control to adjust or correct the operational parameters during disturbances, vector $\mathbf{z}$ is absent, thus, the problem formulation in Eq. (2) reduces to:
$$\psi(\mathbf{x}) = \max_{j \in J}\{g_j(\mathbf{x})\} \text{ s.t. } \mathbf{x} \in \mathbf{T}. \tag{10}$$
It is mathematically equivalent to:
$$\psi(\mathbf{x}) = -\min_{j \in J}\{-g_j(\mathbf{x})\} \text{ s.t. } \mathbf{x} \in \mathbf{T}. \tag{11}$$
Introducing the notation $\varphi_j(\mathbf{x}) = -g_j(\mathbf{x})$, consider for simplicity the case of two constraints $(\varphi_1, \varphi_2)$. From (5) and (9), applying the conjunction R-function we obtain:
$$\psi(\mathbf{x}) = -\min_{j \in J=2}(\varphi_1, \varphi_2) = -R_\wedge(\varphi_1, \varphi_2). \tag{12}$$

We note that the R-conjunction $R_\wedge$ is non-negative when its arguments are non-negative ($\varphi_i \geq 0$). However, R-functions can also handle negative arguments and describe regions in which all arguments are negative, depending on their specific definition and application.

Substituting $x_1, x_2$ with $\varphi_1, \varphi_2$ in (6) and setting the smoothness parameter $\alpha = 1$ to achieve the highest smoothness of the R functions, we obtain:

$$\begin{aligned} R_\wedge(\varphi_1, \varphi_2) &= \varphi_1 \wedge_1 \varphi_2 \equiv \tfrac{1}{2}\left(\varphi_1 + \varphi_2 - \sqrt{\varphi_1^2 + \varphi_2^2 - 2\varphi_1\varphi_2}\right), \\ R_\vee(\varphi_1, \varphi_2) &= \varphi_1 \vee_1 \varphi_2 \equiv \tfrac{1}{2}\left(\varphi_1 + \varphi_2 + \sqrt{\varphi_1^2 + \varphi_2^2 - 2\varphi_1\varphi_2}\right), \end{aligned} \quad (13)$$

which can also be written as

$$\begin{aligned} R_\wedge(\varphi_1, \varphi_2) &= \tfrac{1}{2}(\varphi_1 + \varphi_2 - |\varphi_1 - \varphi_2|), \\ R_\vee(\varphi_1, \varphi_2) &= \tfrac{1}{2}(\varphi_1 + \varphi_2 + |\varphi_1 - \varphi_2|). \end{aligned} \quad (14)$$

It is equivalent to

$$\begin{aligned} R_\wedge(\phi_1, \phi_2) &= \min(\phi_1, \phi_2), \\ R_\vee(\phi_1, \phi_2) &= \max(\phi_1, \phi_2). \end{aligned} \quad (15)$$

Scaling (12) to $J$ constraints, we obtain:

$$\psi(x) = -\min_{j \in J}\{\varphi_j(x)\} = -R_\wedge(\varphi_1, \varphi_2, \dots, \varphi_J) \text{ s.t. } x \in T. \quad (16)$$

Explicitly, the R-conjunction for $J$ constraints is given by:

$$R_\wedge(\varphi_1, \dots, \varphi_d) = \varphi_1(x) \wedge_1 \dots \wedge_1 \varphi_J(x). \quad (17)$$

We conclude that for explicitly given constraints, the optimization problem (10) involving the set of objective functions $\{g_j(x), j = 1, \dots, J\}$ can be solved analytically by applying the R-conjunction function. The solution is expressed as a single explicit function:

$$R_\wedge(x) = R_\wedge(-g_1(x), -g_2(x), \dots, -g_J(x)). \quad (18)$$

From this formulation:

- $R_\wedge(x) \geq 0$ implies that $-g_j(x) \geq 0$ for all $j \in J$, meaning $g_j(x) \leq 0$ for all $j \in J$, ensuring feasibility.
- $R_\wedge(x) < 0$ implies that at least one $-g_j(x) < 0$, meaning at least one $g_j(x) > 0$, indicating infeasibility.

Thus, the feasibility function can be compactly expressed as:

$$\psi(x) = -R_\wedge(x) \text{ s.t. } x \in T. \quad (19)$$

We note the similarity of this problem with the DS identification problem as formulated in Kucherenko et al (2024):

$$DS = DS(x|\, R_\wedge\left(\phi_1(x), \dots, \phi_J(x)\right) \geq 0). \quad (20)$$

4.2. Closed loop scenario

In the closed loop scenario the feasibility function is defined by solving the following optimization problem:

$$\psi(x) = \min_{z}\max_{j\in J}\{g_j(z,x)\} \text{ s.t. } x \in T. \tag{21}$$

Using the methodology outlined in the previous Section the R-conjunction function in this case has the form:

$$R_\wedge(z,x) = R_\wedge(\phi_1,\dots,\phi_J) = \min_{j\in J}\{-g_j(z,x)\}. \tag{22}$$

Hence (21) can be reformulated as

$$\psi(x) = \min_{z}\{-R_\wedge(z,x)\} \text{ s.t. } x \in T. \tag{23}$$

## 5. Test cases.

### 5.1. Open-loop scenario

To demonstrate the performance of the R-feasibility approach, a set of benchmark examples for feasibility analysis is selected from the literature. All examples are low-dimensional, allowing for a graphical representation of the results.

We recall that the standard form of R-functions describe regions in which feasible conditions are expressed using non-negative terms. Hence in this Section we transform constraints $f_i \leq 0$ into $\varphi_i = -f_i \geq 0$ when we apply the R-function method. In all the examined cases, we employ the R-functions methodology by utilizing Eqs. (14) and (17) to determine the feasible region (19): $\psi(x) = -R_\wedge(x)$.

### Example 1. Heat exchanger network problem

The heat exchanger network problem in Eq. (1) from Grossmann et al. (2014) is a minor modification of the seminal example by Saboo and Morari (1984), where the heat capacity flowrate $F_{H1}$ is treated as an uncertain parameter (see Fig. 1). The objective is to assess whether the network is feasible within the range $1 \leq F_{H1} \leq 1.8$ kW/K. By utilizing the relevant heat balance equations, temperature constraints can be formulated in terms of the cooling load $Q_c$, which acts as the control variable or degree of freedom, and the uncertain parameter $F_{H1}$. To ensure feasible operation of the network, the following constraints must be met:

$$\begin{aligned}
f_1 &= -25 + Q_c\left[\frac{1}{F_{H1}} - 0.5\right] + \frac{10}{F_{H1}} \leq 0, \\
f_2 &= -190 + \frac{10}{F_{H1}} + \frac{Q_c}{F_{H1}} \leq 0, \\
f_3 &= -270 + \frac{250}{F_{H1}} + \frac{Q_c}{F_{H1}} \leq 0, \\
f_4 &= 260 - \frac{250}{F_{H1}} - \frac{Q_c}{F_{H1}} \leq 0.
\end{aligned} \tag{18}$$

Grossmann et al. (2014) demonstrated feasibility by plotting the inequalities as a function of the control variable $Q_c$ and the uncertain parameter $F_{H1}$. Their analysis shows that feasible operation is achieved at the extreme points $F_{H1} = 1$ kW/K with $Q_c = 15$ kW and $F_{H1} = 1.8$ kW/K with $Q_c = 227$ kW by adjusting the cooling load. They highlight that without such a plot, one might incorrectly assume the network is feasible across the entire range $1 \leq F_{H1} \leq 1.8$ kW/K. However, as illustrated in Fig. 2 of their paper, at an intermediate value such as $F_{H1} = 1.2$ kW/K with $Q_c = 58.6$ kW, the feasible region is empty. In fact, there exists a non-convex infeasible region between $F_{H1} = 1.118$ and $1.65$ kW/K, with the most significant violation occurring at $F_{H1} = 1.37$ kW/K, which represents the critical point.

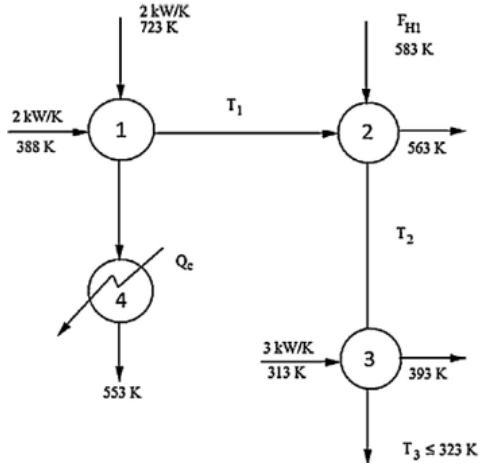

Fig. 1. Heat exchanger network with uncertain heat capacity flowrate, $F_{H1}$ (Grossmann et al., 2014).

By applying the R-function methodology, the feasible region is analytically defined as $\psi(x) = -R_\wedge(x)$, where $x = (F_{H1}, Q_c)$ and

$$R_\wedge(x) = R_\wedge(\phi_1, \phi_2, \phi_3, \phi_4) == \frac{5}{2} - \frac{B}{2} - \frac{A}{2} - \frac{|A - B + 5|}{2},$$

$$A = \frac{4Q_c - 430 F_{H1} + 2F_{H1}|Q_c/2 - 165| - F_{H1} Q_c + 40}{4 F_{H1}},$$

$$B = \frac{|2Q_c - 530 F_{H1} + 500|}{2 F_{H1}}.$$

Curves representing constraints (18) and the feasible region boundary $R_\wedge(\phi_1, \phi_2, \phi_3, \phi_4) = 0$ are shown in Fig. 2. In Fig. 3 a heatmap illustrates the function values within the feasible region $R_\wedge(\phi_1, \phi_2, \phi_3, \phi_4) \geq 0$.

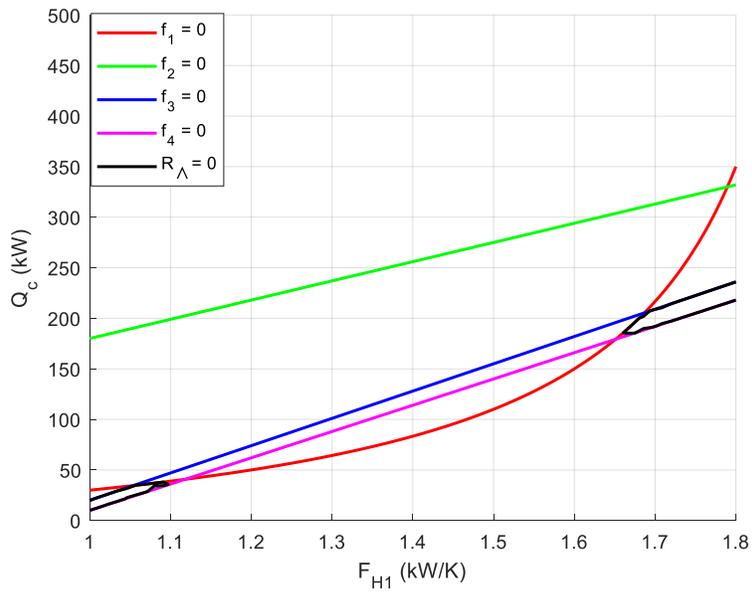

Fig. 2. Curves representing constraints (18) and the feasible region boundary (black curve): $R_\wedge(\phi_1, \phi_2, \phi_3, \phi_4) = 0$. Example 1.

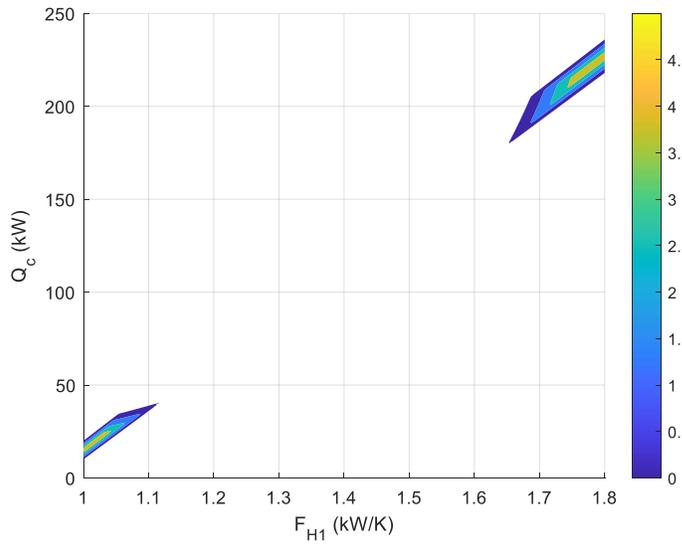

Fig. 3. A heatmap of the feasible region $R_\wedge(\phi_1, \phi_2, \phi_3, \phi_4) \geq 0$. Example 1.

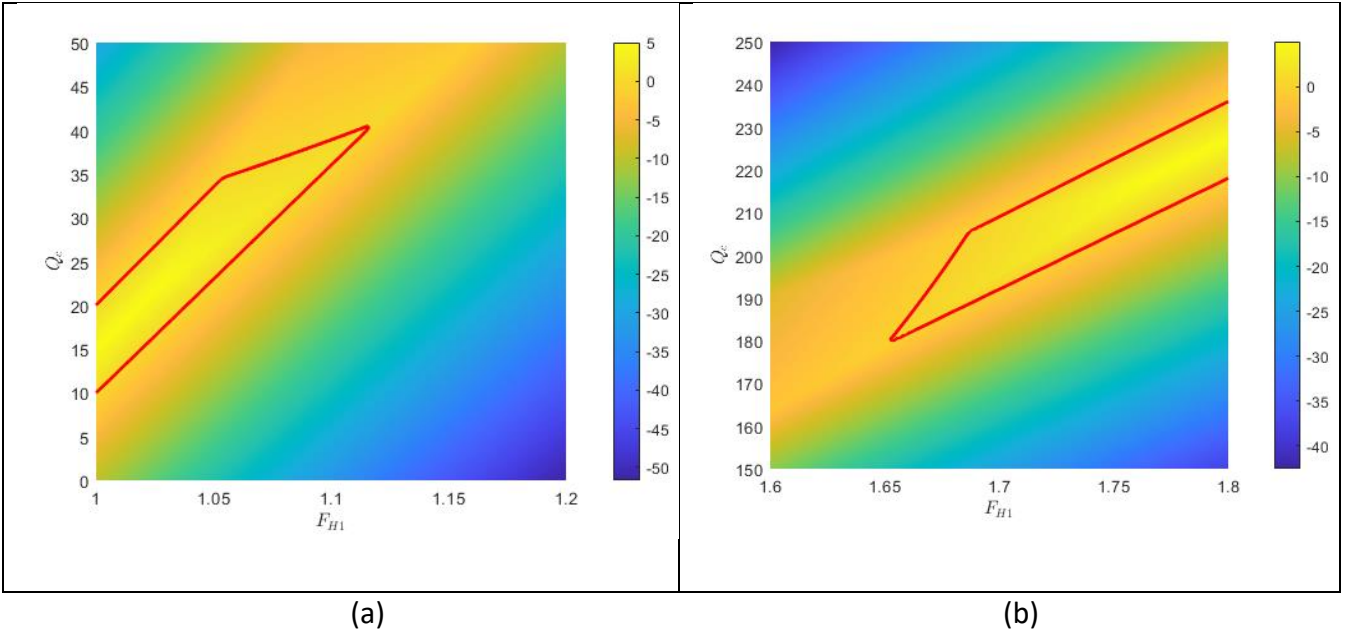

Fig. 4. Enlarged heatmap view of $R_\wedge(\varphi_1, \varphi_2, \varphi_3, \varphi_4)$. (a) Left part of the region; (b) Right part of the region. Red curves indicate the boundary of the feasible region. Example 1.

For operational considerations, if we aim to remove the sharp tip on the right-hand side of the feasible region for $F_{H1} \geq 1.7$, we introduce an additional constraint:
$$f_5 = 1.7 - F_{H1} \leq 0.$$
This modification results in an updated feasible region defined by $R_\wedge(\varphi_1, \varphi_2, \varphi_3, \varphi_4, \varphi_5) \geq 0$ and shown in Fig. 5:

$$R_\wedge(x) = \frac{F_{H1}}{2} - \frac{|A - B + 5|}{4} - \frac{|10F_{H1} + 5|A - B + 5| + 5B + 5A - 42|}{20} - \frac{B}{8} - \frac{A}{16} + \frac{2}{5} \geq 0,$$

$$A = \frac{4Q_c - 430F_{H1} + 2F_{H1}|Q_c/2 - 165| - F_{H1}Q_c + 40}{4F_{H1}},$$

$$B = \frac{|2Q_c - 530F_{H1} + 500|}{2F_{H1}}.$$

Our approach provides a more detailed and explicitly formulated solution compared to the method presented by Grossmann et al. (2014), which is primarily descriptive.

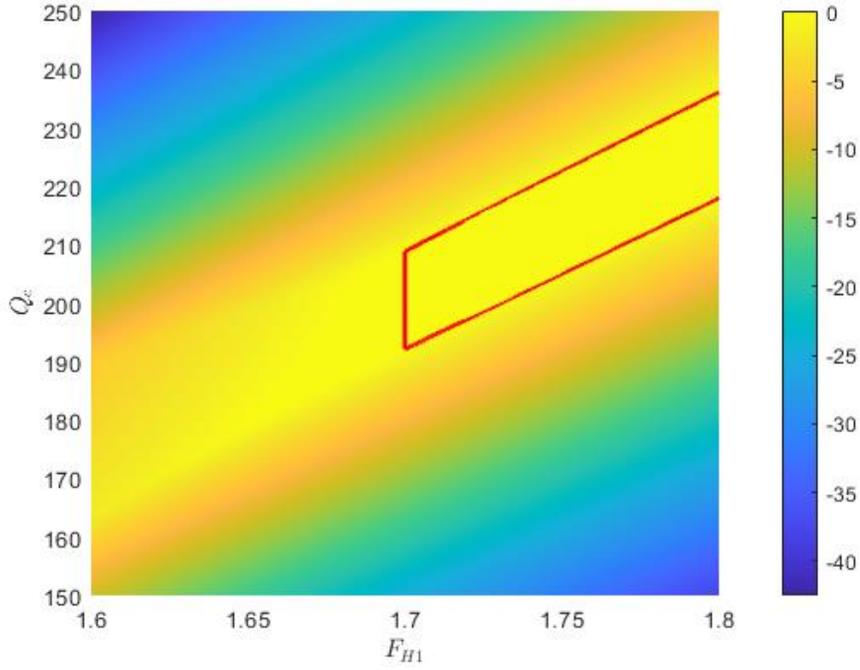

Fig. 5. Enlarged view of a heatmap displaying the right part of the trimmed feasible region. The red curve represents the boundary of the feasible region. Example 1.

Example 2. Three nonlinear convex constraints

This example was initially proposed by Ierapetritou (2001) and subsequently revisited by Goyal and Ierapetritou (2002). It involves two uncertain parameters and three nonlinear convex constraints. The feasibility problem in this case is defined by the following set of inequalities:

$$f_1 = \theta_2 + \theta_1^2 - \theta_1 - 40 \leq 0,$$
$$f_2 = \theta_1^2 + \theta_1 - \theta_2 - 2 \leq 0, \qquad (19)$$
$$f_3 = \theta_2 - 4 \times \theta_1 - 30 \leq 0.$$

The authors proposed a method for identifying operating envelopes by iteratively refining the feasible region's boundary approximation using the simplicial approximation technique (Director and Hachtel, 1977). This method identifies boundary points, which are then used to construct a convex hull inscribed within the feasible region, serving as a lower bound of the feasible space. The outer boundary is determined by generating tangent planes at these points, providing an upper bound of the feasible region.

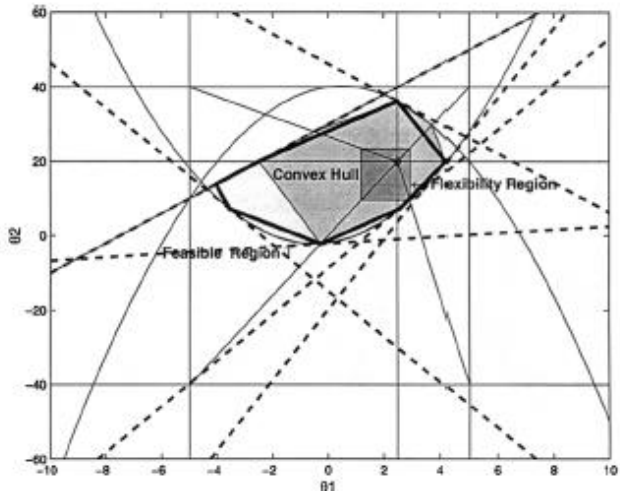 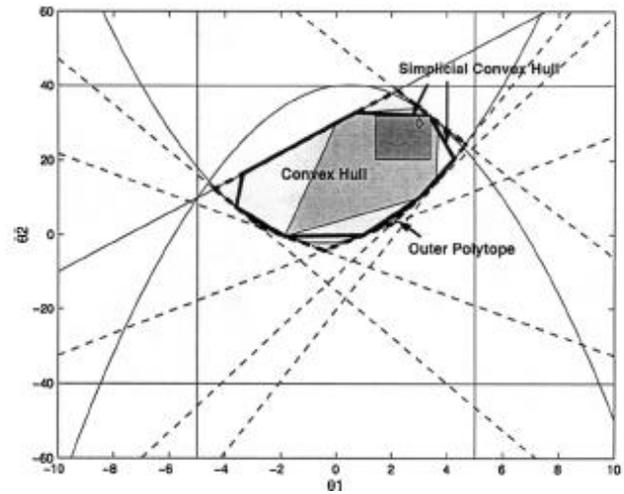

(a) (b)

Figure 6. (a) Simplicial convex hull with nominal points of $(\theta 1, \theta 2) = (2.5, 20)$; Simplicial convex hull with nominal points of $(\theta 1, \theta 2) = (3, 30)$ (b). Example 2 (Goyal and Ierapetritou, 2002).

Solving the feasibility problem using simplicial approximation involved solving seven nonlinear and five linear programs across all iterations, achieving convergence after four simplicial iterations. On average, the NLPs required 25 iterations per run, while the LPs converged within five iterations (Goyal and Ierapetritou, 2002).

To demonstrate the independence of the proposed technique from nominal point location, the nominal points were shifted from $(\theta_1, \theta_2) = (2.5, 20)$ to $(3, 30)$, and the analysis was repeated. The simplicial approximation approach resulted in a slight increase in the inscribed convex hull volume to 178.56 units and a slight decrease in the outer convex polytope volume to 204.83 units (Fig 6).

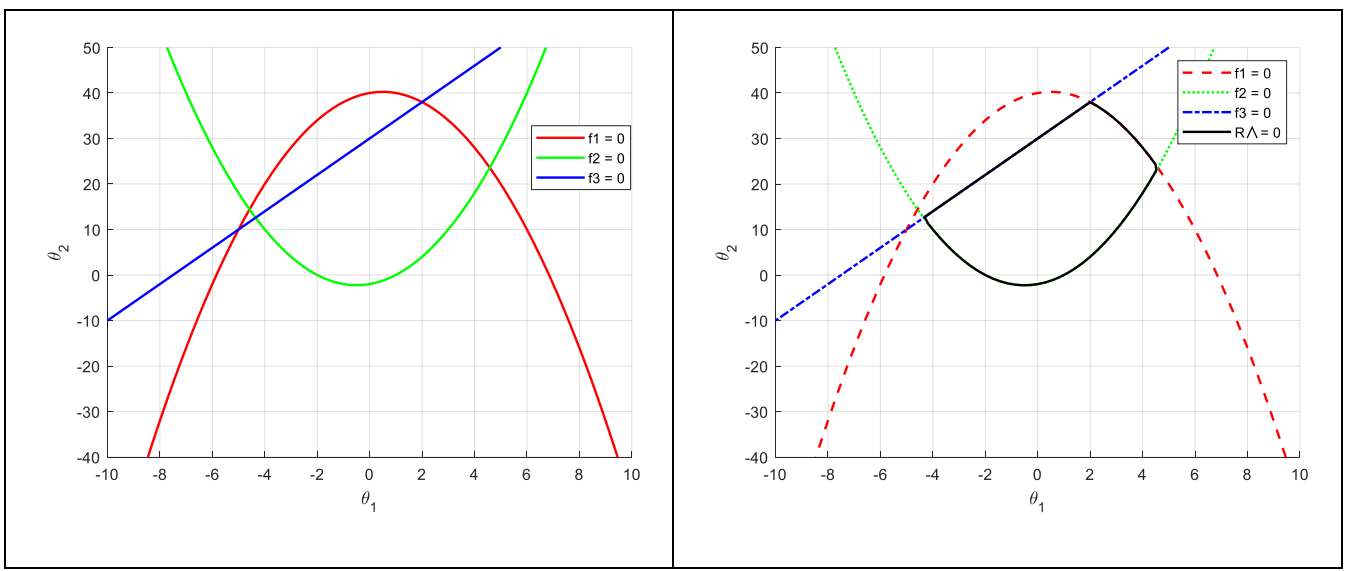

(a) (b)

Fig. 7. Curves representing constraints (19) and the feasible region boundary (black curve): $R_\wedge(\phi_1, \phi_2, \phi_3) = 0$. Example 2.

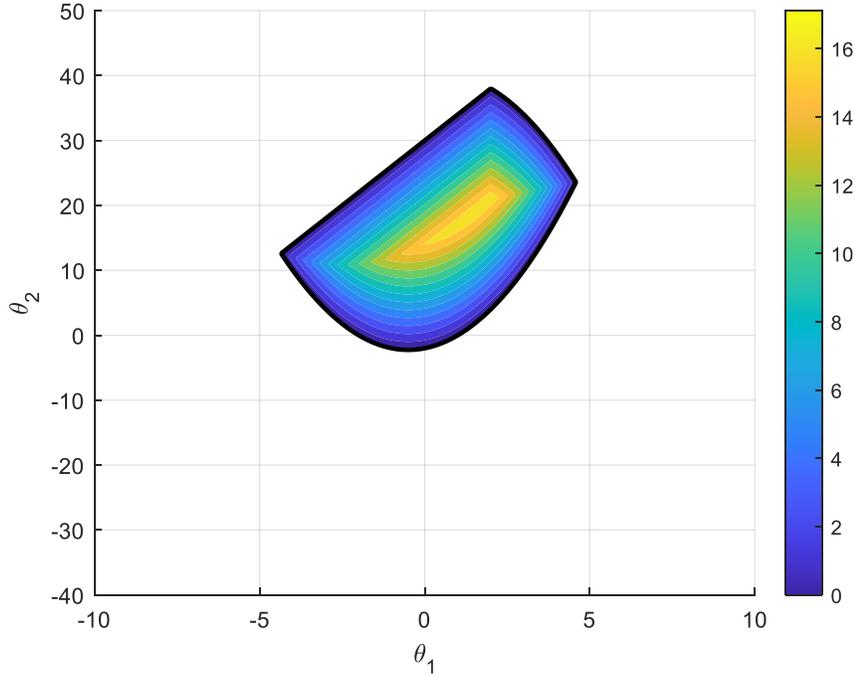

Fig. 8. A heatmap of $R_\wedge(\phi_1, \phi_2, \phi_3) \geq 0$. Example 2.

Using R-functions the feasible region, expressed with a negative sign as $\psi(x) = -R_\wedge(x)$, with $x = (\theta_1, \theta_2)$ is analytically defined as:

$$R_\wedge(\phi_1, \phi_2, \phi_3) = 2\theta_1 - \frac{\theta_2}{2} - \frac{|2\theta_1 - 2\theta_2 + 38|}{4} - \frac{\left|4\theta_1 - \theta_2 + \frac{|2\theta_1 - 2\theta_2 + 38|}{2} + \theta_1^2 + 9\right|}{2} - \frac{\theta_1^2}{2} + \frac{51}{2} \geq 0.$$

The curves depicting the constraints and the boundary of the feasible region $R_\wedge(\phi_1, \phi_2, \phi_3) = 0$ are displayed in Fig. 7. Figure 8 presents a heatmap showing the function values within the feasible region $R_\wedge(\phi_1, \phi_2, \phi_3) \geq 0$. The calculated volume of this true feasible region is 199.4 units.

Example 3. The convex hull approach in the presence of a nonconvex constraint

Banerjee and Ierapetritou (2005) proposed the following set of convex and nonconvex constraints to characterize a feasible region:

$$f_1 = \theta_2 - 2\theta_1 - 15 \leq 0,$$
$$f_2 = \frac{\theta_1^2}{2} + 4\theta_1 - 5 - \theta_2 \leq 0, \quad (20)$$
$$f_3 = \theta_2(6 + \theta_1) - 80 \leq 0,$$

where $\theta_1$ and $\theta_2$ are the uncertain parameters in the range $-10 \leq \theta_1 \leq 10$ and $-15 \leq \theta_2 \leq 15$. They observed that the convex hull approach is restricted to convex and one-dimensional quasi convex feasible regions, with its effectiveness declining when nonconvex constraints are present. To demonstrate this limitation, an additional nonconvex constraint was introduced:

$$f_4 = 10 - \frac{(\theta_1 - 4)^2}{5} - 2\theta_2^2 \leq 0.$$

To determine the range of parameters within which a given process remains viable requires obtaining a mathematical representation of the feasible region within the parameter space, which is defined by the process constraints. Banerjee and Ierapetritou (2005) suggested to conceptualize the feasible region as a geometric object, with its shape or boundary estimated using the $\alpha$-shape technique. A possible approach to surface reconstruction involves using an input set of points that represent the object whose surface needs to be determined. However, the R-function-based methodology eliminates this step by directly providing an explicit mathematical representation of the feasible region within the parameter space, as defined by the process constraints.

In the $\alpha$-shape method, after generating a representative set of data points for the feasible region, the $\alpha$-shape is constructed based on an $\alpha$ estimate obtained from the minimum spanning tree of the data set. The $\alpha$-shape identifies the points in the input data set that define the object's surface. These points are then connected-by line segments in two dimensions or by triangles in three dimensions-forming a polygon that encloses the feasible region. Generating points that capture the shape of the feasible region typically involves uniform sampling across the entire parameter space, without distinguishing between feasible and infeasible regions. To address this limitation, Banerjee and Ierapetritou (2005) proposed an alternative approach that reformulates the sampling process as an optimization problem, which is then solved using a genetic algorithm (GA). For the considered test case Banerjee and Ierapetritou (2005) observed that 1) only 3,064 out of 40,000 function calls were unique, resulting in 938 feasible points; 2) the value of $\alpha$ plays a critical role in determining the level of detail captured by the $\alpha$-shape. For the 938 sampled points, the computed $\alpha$ value was 25, which provided a sufficiently accurate representation of the object's nonconvex nature (Fig. 9a). However, if the value of $\alpha$ is increased, the $\alpha$-shape prediction tends to lose its accuracy (Fig. 9b). The authors noted that applying the proposed technique in higher dimensions poses significant computational challenges.

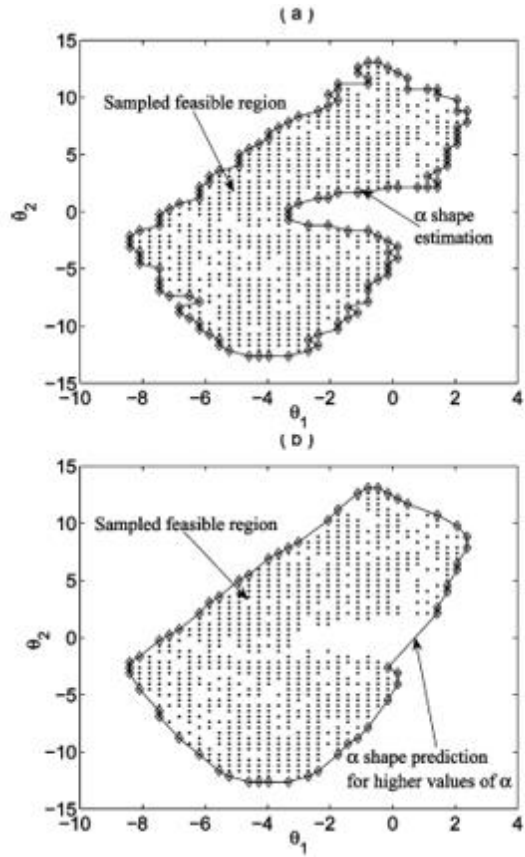

Figure 9. Performance of the $\alpha$ shape method in predicting the feasible space using different $\alpha$ values: (a) $\alpha = 25$; (b) $\alpha = 1000$. (Banerjee and Ierapetritou, 2005). Example 3.

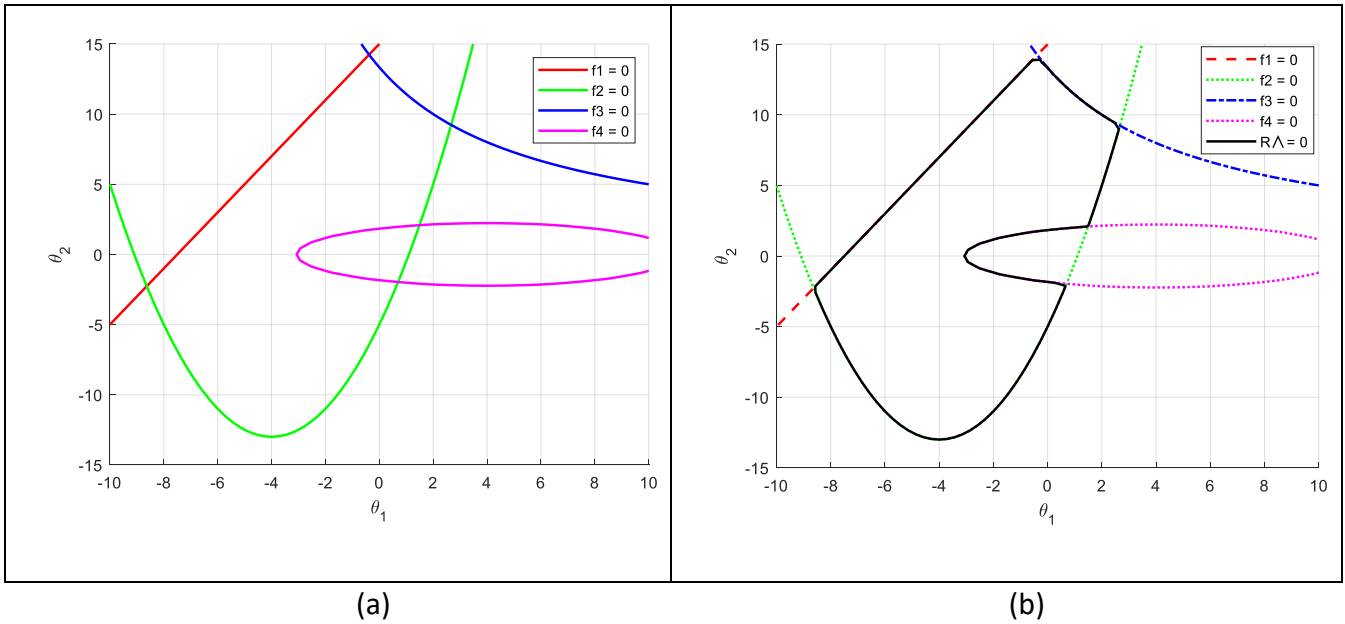

Fig. 10. Curves representing constraints (20) and the feasible region boundary (black): $R_\wedge(\phi_1, \phi_2, \phi_3, \phi_4) = 0$. Example 3.

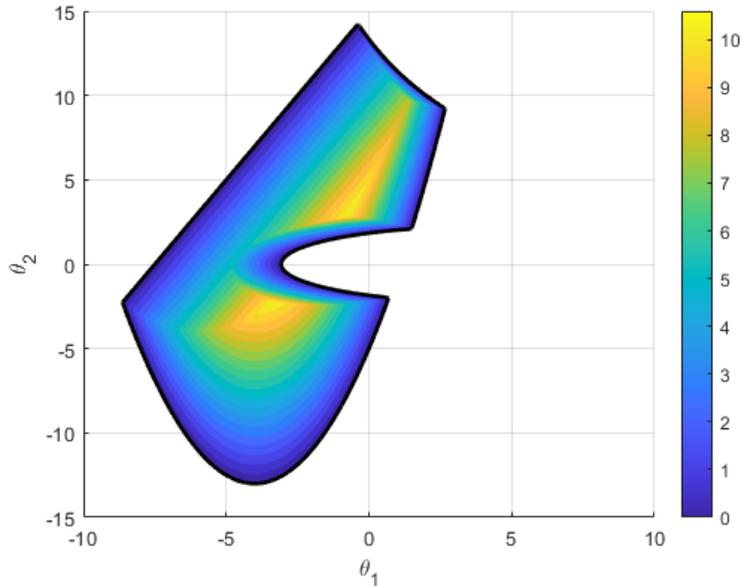

Fig. 11. A heatmap of the feasible region $R_\wedge(\phi_1, \phi_2, \phi_3, \phi_4) \geq 0$. Example 1.

The feasible region, represented with a negative sign as $\psi(x) = -R_\wedge(x)$, where $x = (\theta_1, \theta_2)$, is analytically defined as:

$R_\wedge(\phi_1, \phi_2, \phi_3, \phi_4) =$

$$= \frac{-15\theta_1^2 - 36\theta_1 + 10\theta_2^2 - 60\theta_2 - 10\theta_1\theta_2 - 5 \cdot \left|6\theta_1 - 2\theta_2 + \frac{\theta_1^2}{2} + 10\right| - 5 \cdot \left|\frac{(\theta_1 - 4)^2}{5} + \theta_2(\theta_1 + 6) + 2\theta_2^2 - 90\right| + 466}{20}$$

The curves depicting the constraints and the boundary of the feasible region $R_\wedge(\phi_1, \phi_2, \phi_3, \phi_4) = 0$ are illustrated in Fig. 10. Figure 11 presents a heatmap displaying the $R$-function values within the feasible region $R_\wedge(\phi_1, \phi_2, \phi_3) \geq 0$. The volume of the feasible region, calculated using the $R$-function method, is 152.17 units. Banerjee and Ierapetritou (2005) did not provide the volume obtained through their approach.

Table. 1. Comparison of the $\alpha$-Shape and R-Function Methods

| Aspect | $\alpha$-Shape Method | R-Function Method |
| --- | --- | --- |
| **Step 1: Sampling the Feasible Region** | Sample the feasible space to obtain a set of points representing the region. | No need to sample. |
| **Step 2: Constructing the Feasible Region** | Build an $\alpha$-shape using the sampled points to identify the boundary points. | Construct an analytical $R_\wedge$ function from the given constraints. |

| Step 3: Creating the Feasible Region Boundary | Connect the boundary points to form a polygonal representation of the feasible region. | Solve the equation $R_\wedge = 0$. |
| Step 4: Feasibility Check | Perform a point-in-polygon test to check if a point is inside the constructed polygon. | Compute the function $R_\wedge(x)$ to determine its sign (feasibility) of point $x$. |

Table 1 highlights the main differences between compared methods, showing that the $\alpha$ shape method relies on sampled points and polygon construction, while the R-function method directly leverages analytical expressions to define and evaluate feasibility efficiently. The feasibility check using the $\alpha$-shape method also involves a computationally expensive point-inpolygon test to determine whether a point lies within the constructed polygon. In contrast, the R-function method requires only a simple function evaluation check $R_\wedge(x) \geq 0$ to assess the feasibility of a point $x$.

Example 4. Disjoint feasible region

This example from Boukouvala and Ierapetritou (2012) considered a situation where the feasible region consists of disjoint components:

$$\begin{aligned}
f_1 &= 4\theta_1^2 - 2.1\theta_1^4 + \frac{\theta_1^6}{3} + \theta_1\theta_2 - 4\theta_2^2 + 4\theta_2^4 \leq 0, \\
f_2 &= 2\theta_1 - \theta_2 - 3 \leq 0, \\
f_3 &= -0.8\theta_1 + \theta_2 - 1.8 \leq 0, \\
\theta_1 &\in [-2,2], \\
\theta_2 &\in [-1,1].
\end{aligned} \quad (21)$$

The first constraint is the well-known six-hump camel test function, commonly used as a benchmark for optimization methods. It generates a multimodal surface, enabling the simulation of a disjoint feasible region.

Boukouvala and Ierapetritou (2012) employed a uniform sampling strategy for the initial experimental design. After two iterations and a total of 29 design samples, the feasible region predicted using a Kriging surrogate response surface is presented in Fig. 12. While the predicted feasible space is underestimated, the proposed approach allows to handle discontinuities and non-linearity within the feasible region.

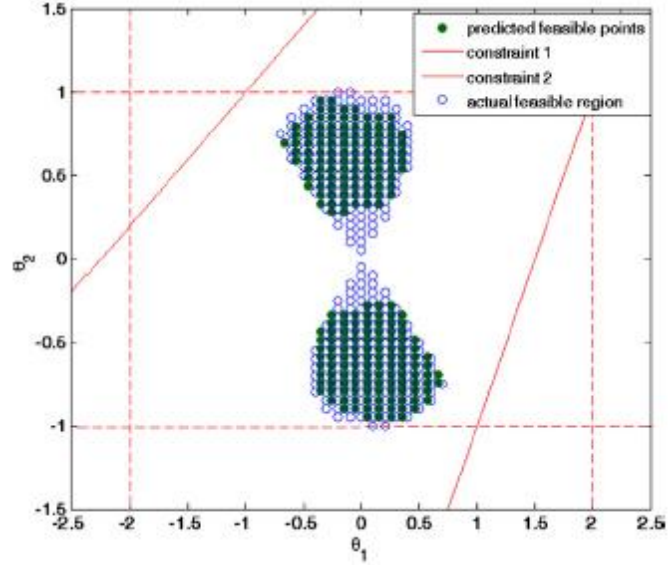

Fig. 12. Predicted vs. real feasible region of a non-linear disjoint problem (Boukouvala and Ierapetritou, 2012). Example 4.

The feasible region, expressed with a negative sign as $\psi(x) = -R_\wedge(x)$, with $x = (\theta_1, \theta_2)$ is analytically defined as:

$$R_\wedge(\phi_1, \phi_2, \phi_3) = \theta_2^2 - \frac{\theta_2}{4} - \frac{1}{2}\left|\frac{9\theta_1}{5} - \frac{3\theta_2}{2} + \frac{A}{2} + \frac{\theta_1\theta_2}{2} + 2\theta_1^2 - 2\theta_2^2 - \frac{21\theta_1^4}{20} + 2\theta_2^4 + \frac{\theta_1^6}{6} + \frac{3}{10}\right|$$

$$- \frac{A}{4} - \frac{\theta_1\theta_2}{4} - \theta_1^2 - \frac{\theta_1}{10} + \frac{21\theta_1^4}{40} - \theta_2^4 - \frac{\theta_1^6}{12} + \frac{33}{20},$$

$$A = \left|\theta_2 - 2\theta_1 + \theta_1\theta_2 + 4\theta_1^2 - 4\theta_2^2 - \frac{21\theta_1^4}{10} + 4\theta_2^4 + \frac{\theta_1^6}{3} + 3\right|.$$

Curves representing constraints (21) and the feasible region boundary $R_\wedge(\phi_1, \phi_2, \phi_3) = 0$ are shown in Fig. 13. In Fig. 14 a heatmap illustrates the function values within the feasible region $R_\wedge(\phi_1, \phi_2, \phi_3) \geq 0$. We observe that the feasible region is precisely defined, unlike in the method proposed by Boukouvala and Ierapetritou (2012). The computed volume of the feasible region is 1.46 units. Boukouvala and Ierapetritou (2012) did not report the volume obtained using their method.

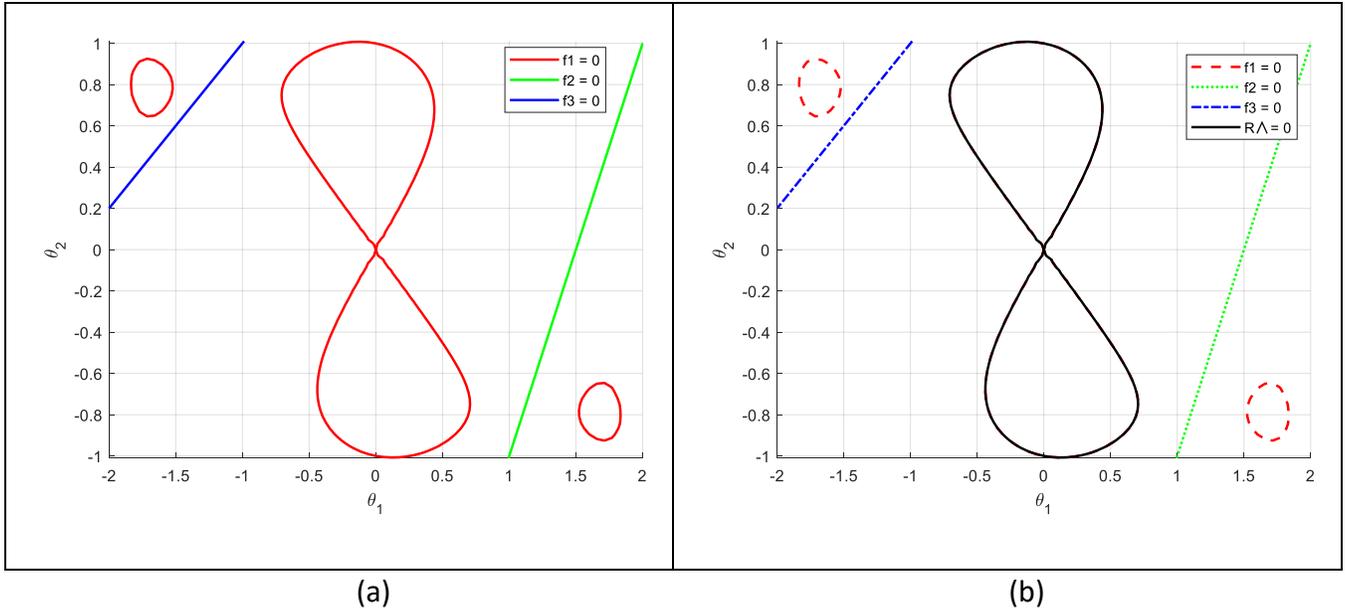

Fig. 13. Curves representing constraints (21) (a); the feasible region boundary (black): $R_\wedge(\phi_1, \phi_2, \phi_3) = 0$. (b) Example 4.

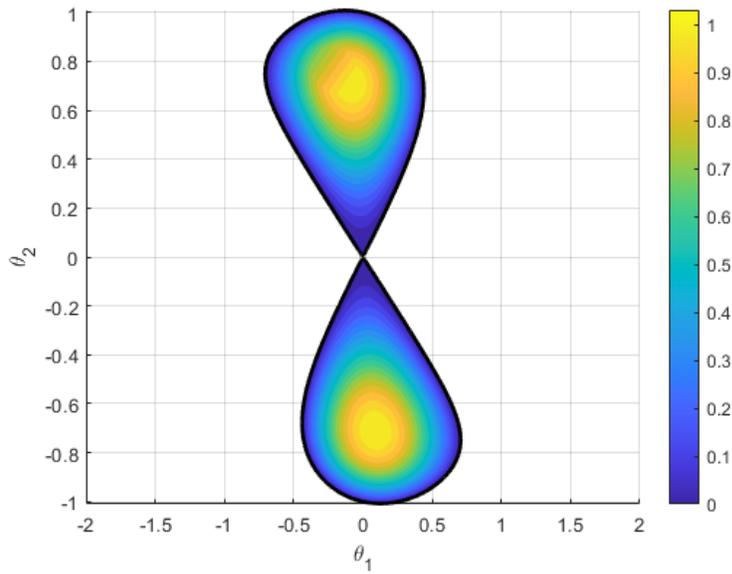

Fig. 14. A heatmap of $R_\wedge(\varphi_1, \varphi_2, \varphi_3) \geq 0$. Example 4.

## 5.2. Closed-loop scenario

In this section, the problems involve a single control variable, $z$, along with process parameters.

Example 5.

Consider the following set of two constraints which involve one control variable $z$ and one parameter $\theta$ from Halemane and Grossmann (1983):

$$f_1 = -z + \theta \leq 0,$$
$$f_2 = z - 2\theta + 2 - d \leq 0, \quad (22)$$
$$1 \leq \theta \leq 2.$$

The feasibility function $\psi(d, \theta)$ is obtained by solving the following problem for $d = 0.5$ and $1 \leq \theta \leq 2$:

$$\psi(d, \theta) = \min_z u$$
$$\text{s.t. } u \geq f_1 = -z + \theta$$
$$u \geq f_2 = z - 2\theta + 2 - d.$$

Fig. 1a in Halemane and Grossmann, 1983 illustrates the feasible region in the $(z, \theta)$ space for a design where $d = 0.5$. As shown in the figure, the feasible region expands as $\theta$ increases: $\theta = 1$ is infeasible, $\theta = 1.5$ has feasibility at a single point of $z$, and $\theta = 2$ exhibits a finite feasible region.

From the results for the feasibility function $\psi(d, \theta)$ illustrated in Fig. 1b (Halemane and Grossmann, 1983) it follows that $\psi = 0$ at $\theta = 1.5$, indicating a single feasibility point. Negative $\psi$ values signify finite feasible regions, such as at $\theta = 2$, while positive values indicate infeasibility, as seen at $\theta = 1$, where $\psi$ reaches its maximum, marking the critical point. $\psi$ decreases monotonically as $\theta$ increases, representing the expansion of the feasible region. The quantity $\psi$ serves as a measure of the feasible region's extent, reflecting its operability. This region corresponds to the projection of the overall feasible domain in the $(z, \theta, d)$ space onto the $z$-space, with fixed $d$ and $\theta$ values.

The R-function $R(z, \theta, d)$ that describes the feasible region in the $(z, \theta, d)$ space has the form:

$$R(z, \theta, d) = R_\wedge(\phi_1, \phi_2) = \frac{d}{2} + \frac{\theta}{2} - \frac{|d + 3\theta - 2z - 2|}{2} - 1.$$

$R(z, \theta, d) \geq 0$ defines the feasible region, and $R(z, \theta, d) < 0$ defines the infeasible region. For a fixed value of $d = 0.5$, the following R-function describes the feasible region in the $(z, \theta)$ space (Fig. 15a):

$$R_\wedge(\phi_1, \phi_2) = \frac{\theta}{2} - \frac{|4z - 6\theta + 3|}{4} - \frac{3}{4}.$$

The feasibility function $\psi(d, \theta)$ is obtained by solving the following problem:

$$\psi(\theta) = \min_z(-R_\wedge(z, \theta)),$$
$$1 \leq \theta \leq 2$$

with $R(z, \theta) = R_\wedge(\phi_1, \phi_2)$. The solution to this problem is found explicitly:

$$\psi(\theta) = \frac{3}{4} - \frac{\theta}{2},$$

with the optimal and feasible control $z^* = \frac{6\theta - 3}{4}$. This is illustrated in Fig. 15b. Fig. 16 presents a heatmap of $R_\wedge(\varphi_1, \varphi_2)$ along with the curve for the boundary $R_\wedge(\phi_1, \phi_2) = 0$.

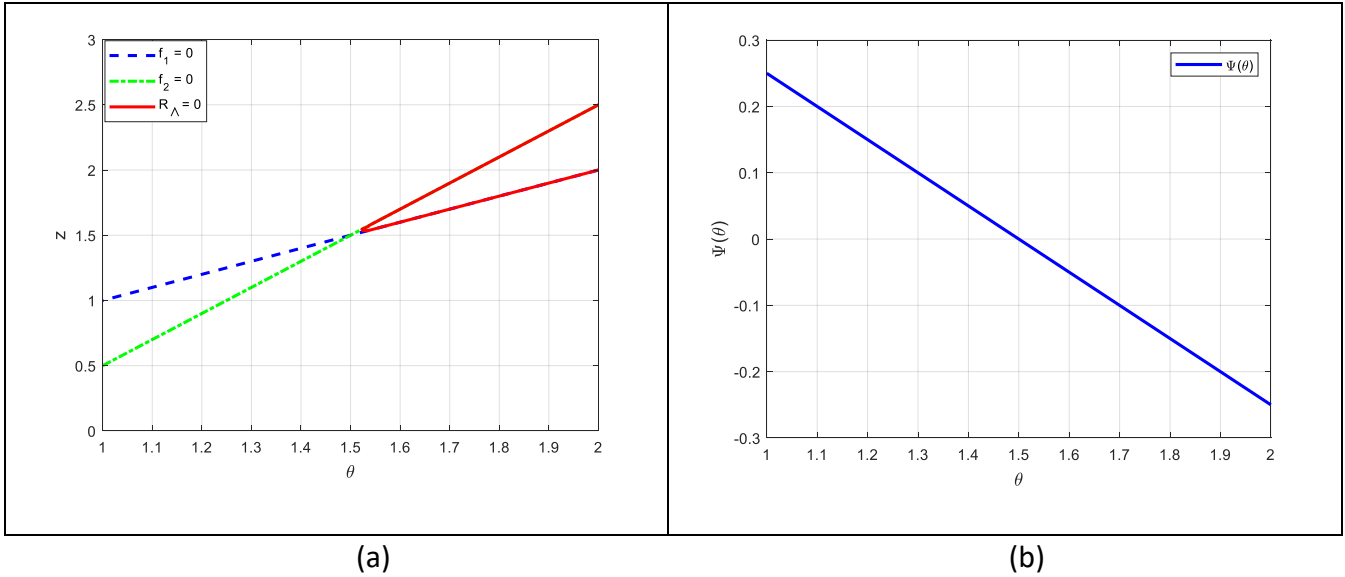

(a)           (b)

Fig. 15. Curves representing constraints (22) and the boundary $R_\wedge(\phi_1, \phi_2) = 0$ (a); the feasibility function $\psi(\theta)$ (b). Example 5.

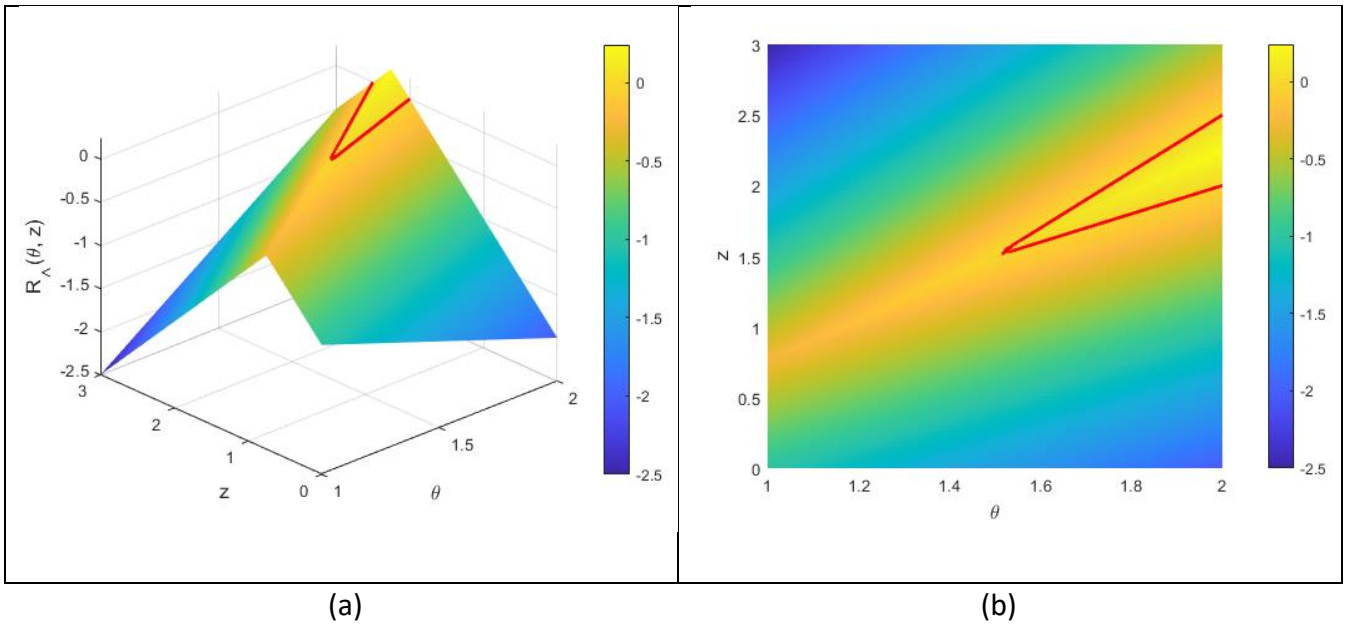

(a)           (b)

Fig. 16. A heatmap of $R_\wedge(\phi_1, \phi_2)$ (a). A heatmap of 2D projection $R_\wedge(\phi_1, \phi_2)$. Red curves represent the boundary $R_\wedge(\phi_1, \phi_2) = 0$ (b). Example 5.

Example 6.

This example from Halemane and Grossmann (1983) illustrates that introducing a third constraint to Example 5 leads to multiple critical $\theta$-points:

$$\begin{aligned} f_1 &= -z + \theta \leq 0, \\ f_2 &= z - 2\theta + 2 - d \leq 0, \\ f_3 &= -z + 6\theta - 9d \leq 0, \\ 1 &\leq \theta \leq 2. \end{aligned} \quad (23)$$

The feasible region for this set of constraints is illustrated in Fig. 17a (Halemane and Grossmann, 1983) for $d = 1$, with the corresponding function $\psi$ shown in Fig. 17b. Notably, $\psi$ is nondifferentiable at $\theta = 9/5$ and exhibits two local maxima at $\theta = 1$ and $\theta = 2$. Fig. 17a shows that the feasible region contracts at the boundary points $\theta = 1$ and $\theta = 2$, while expanding toward the interior point $\theta = 9/5$. Consequently, two critical points must be considered for design, corresponding to the extreme values within the given range $1 \leq \theta \leq 2$.

For a fixed value of $d = 1.0$ the following R-function $R_\wedge(z, \theta)$ describes the feasible region in the $z$-$\theta$ space (Fig. 17a):

$$\begin{aligned} R_\wedge(\phi_1, \phi_2, \phi_3) &= \frac{z}{2} - \frac{3\theta}{2} - \frac{A}{4} - \frac{B}{2} + 2, \\ A &= |2z - 8\theta + 10|, \\ B &= \left|\theta + z + \frac{A}{2} - 4\right|. \end{aligned}$$

Solution to the minimization problem

$$\psi(\theta) = \min_{z} (-R_\wedge(z, \theta)),$$
$$1 \leq \theta \leq 2$$

is found analytically:

$$\psi(\theta) = \frac{\theta}{2} - \frac{1}{2} - \frac{|5\theta - 9|}{2}$$

(Fig. 17b). The optimal and feasible control $z^* = 4\theta - 5$. Fig. 18 presents a heatmap of $R_\wedge(\varphi_1, \varphi_2, \phi_3)$ along with the curve for the boundary $R_\wedge(\phi_1, \phi_2, \phi_3) = 0$.

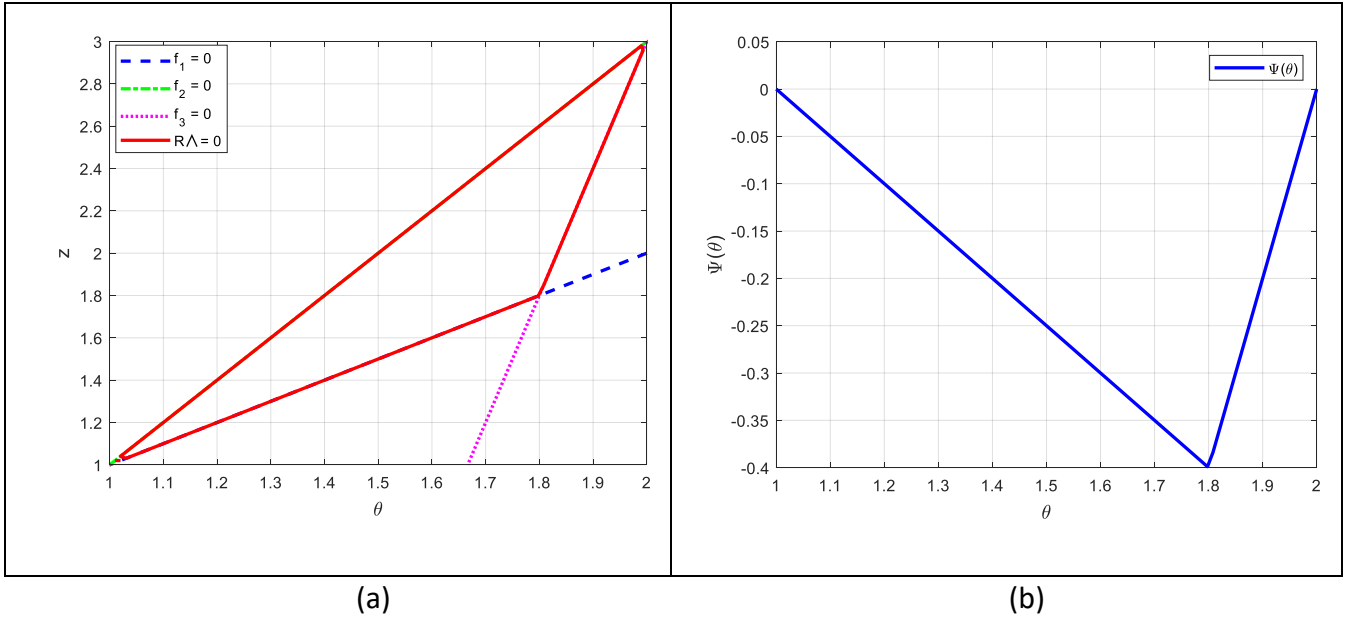

(a)                                           (b)

Fig. 17. Curves representing constraints (23) and the boundary $R_\wedge(\phi_1, \phi_2, \phi_3) = 0$ (a); the feasibility function $\psi(\theta)$ (b). Example 6.

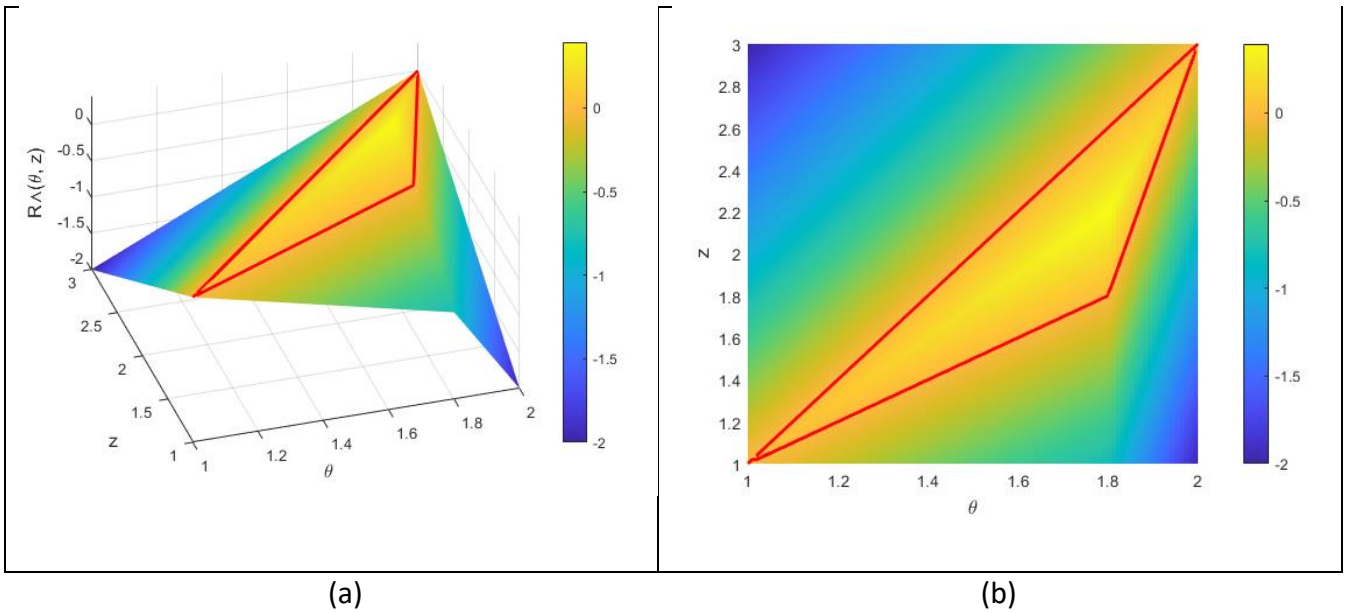

(a)                                           (b)

Fig. 18. A heatmap of $R_\wedge(\phi_1, \phi_2, \phi_3)$ (a). A heatmap of 2D projection $R_\wedge(\phi_1, \phi_2, \phi_3)$. Red curves represent the boundary $R_\wedge(\phi_1, \phi_2, \phi_3) = 0$ (b). Example 6.

Example 7.

The third test case in this Section is taken from Ierapetritou (2001). The constraints describing the feasible region of the design $(d_1, d_2)$ with three uncertain parameters have the following form:

$$f_1 = -z - \theta_1 + 0.5\theta_2^2 + 2.0\theta_3^2 + d_1 - 3d_2 - 8 \leq 0,$$
$$f_2 = -z - \frac{\theta_1}{3} - \theta_2 - \frac{\theta_3}{3} + d_2 + \frac{8}{3} \leq 0, \qquad (24)$$
$$f_3 = z + \theta_1\theta_1 - \theta_2 - d_1 + \theta_3 - 4 \leq 0.$$

Here $z$ is the control variable, and $\theta_1, \theta_2, \theta_3$ are uncertain parameters with a nominal value of $\theta_1^N = \theta_2^N = \theta_3^N = 2$ and expected deviations of $\Delta\theta_1^\pm = \Delta\theta_2^\pm = \Delta\theta_3^\pm = 2$. The design considered further corresponds to $(d_1, d_2) = (3,1)$. In Ierapetritou (2001) the feasible region was approximated by a convex hull with a volume of 12.16 units. Goyal and Ierapetritou (2002) implemented a method based on simplicial approximation, where the feasible region was approximated by a convex hull, resulting in an increased volume of 14.69 units. This approximation is depicted in Fig. 19 (corresponds to Fig 5. in Goyal and Ierapetritou, 2002) as a polyhedron with darker edges. The outer convex polytope, with a volume of 25.29 units, is represented in Fig. 19 with dark dotted edges. This approach required solving 14 nonlinear programs (NLPs) and 44 linear programs (LPs), achieving convergence in ten simplicial iterations. On average, the NLPs required ten iterations per run, while the LPs converged in seven.

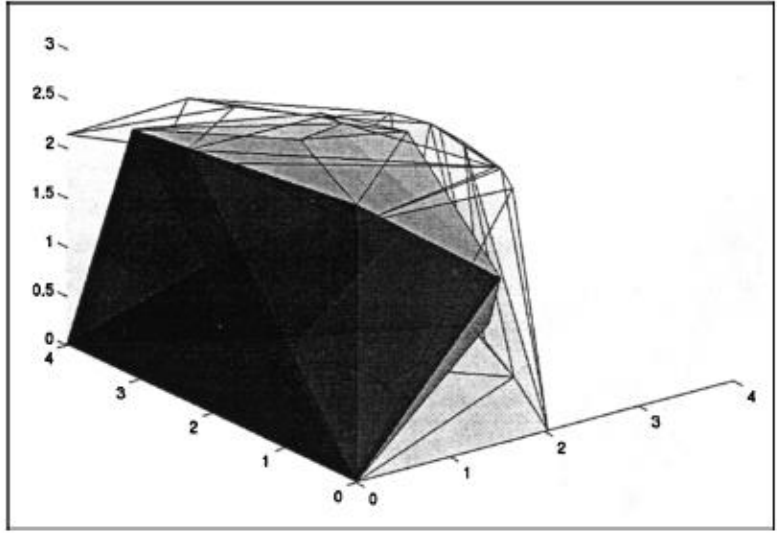

Fig. 19. Simplicial convex hull for a 3-D feasible region (Goyal and Ierapetritou, 2002).

For fixed values $(d_1, d_2) = (3,1)$ the following R-function describes the feasible region in the $z$-$\theta$ space:

$$R(\theta_1, \theta_2, \theta_3, z) == \frac{7\theta_1}{12} + \frac{\theta_2}{2} - \frac{\theta_3}{6} + \frac{z}{2} - \frac{|A|}{4} - \frac{|B|}{2} - \frac{\theta_1^2}{4} - \frac{\theta_2^2}{4} - \theta_3^2 + \frac{29}{6},$$
$$A = \theta_1^2 + \frac{\theta_1}{3} + \frac{4\theta_3}{3} + 2z - \frac{32}{3},$$
$$B = \frac{5\theta_1}{6} - \theta_2 + \frac{\theta_3}{3} + z + \frac{|A|}{2} + \frac{\theta_1^2}{2} - \frac{\theta_2^2}{2} - 2\theta_3^2 + \frac{19}{3}.$$

The solution to the optimization problem

$$\psi(\theta_1, \theta_2, \theta_3) = \min_z \left(-R(\theta_1, \theta_2, \theta_3, z)\right),$$

$$0 \leq \theta_i \leq 4, \, i=1,2,3$$

is found analytically:

$$\psi(\theta_1, \theta_2, \theta_3) = \frac{\theta_1}{2} + \frac{\theta_2}{2} - \frac{\theta_3}{2} - \frac{|3\theta_2^2 + 6\theta_2 + 12\theta_3^2 + 2\theta_3 - 4\theta_1 - 70|}{12} - \frac{\theta_1^2}{2} - \frac{\theta_2^2}{4} - \theta_3^2 + \frac{15}{2}.$$

The optimal and feasible control $z^*$ is reached at $z^* = \frac{32 - 3\theta_1^2 - \theta_1 - 4\theta_3}{6}$. The surface of the feasibility function shown in Fig. 20 resembles the outer convex polytope in Fig. 19 but provides a significantly more accurate solution to the minimization problem above.

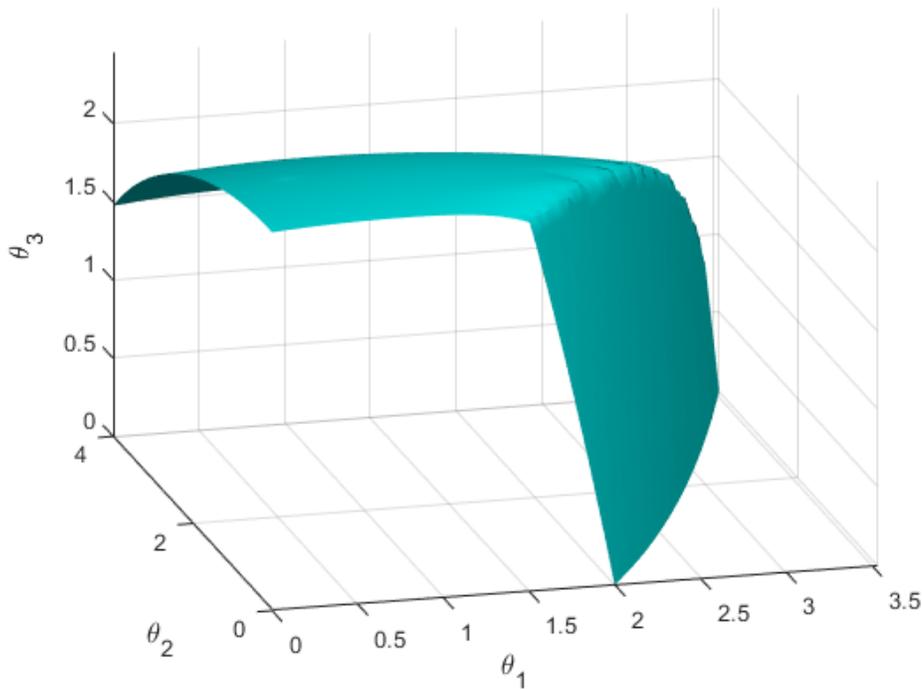

Fig. 20. Surface of the feasibility function $\psi(\theta_1, \theta_2, \theta_3)=0$. Example 6.

The volume of the feasible region is 23.83 units, which is close to the 25.29-unit volume of the outer convex polytope computed by Goyal and Ierapetritou (2002).

## 6. Advantages of using the R-function methodology

The R-function methodology offers several distinct advantages over other existing approaches:
1. It represents the feasibility region using a single implicit algebraic expression.
2. It can handle multiple disjoint feasible regions of any complexity.
3. It seamlessly scales to high-dimensional problems.
4. It is not limited to constraints defined by polynomial functions.
5. It identifies the feasibility region without requiring sampling or optimization steps, thereby avoiding numerical and approximation errors.
6. It is computationally efficient, easy to implement, and adaptable to both symbolic and numerical computation, making it well-suited for high-dimensional applications.

The R-function methodology can leverage symbolic computation to simplify and automate the derivation of implicit functions and the construction of complex geometric or algebraic representations. Software tools like MATLAB, Maple, Mathematica, and Python (e.g., SymPy) are particularly useful for this purpose, as they allow for the exact manipulation of mathematical expressions without relying on numerical approximations. However, when symbolic computation tools are unavailable, the R-function methodology can still be implemented using numerical methods such as grid-based sampling or Monte Carlo techniques to approximate implicit functions and their combinations. Although these numerical approaches may be less precise, they provide flexibility and can be applied in environments where symbolic computation is impractical.

## 6.1. Exploration of the Identified feasible regions

The R-function methodology simplifies feasibility testing by requiring only a straightforward check of whether $R(\varphi_1, \ldots, \varphi_m) \geq 0$ to determine if a point lies within the feasible region. In contrast, traditional approaches approximate the feasible region with a polygon, relying on computationally expensive point-in-polygon tests. One such method is the ray tracing algorithm proposed by Banerjee and Ierapetritou (2005). This algorithm determines whether a point is inside by counting the intersections of a semi-infinite ray drawn from the point with the polygon's edges. An odd number of intersections indicates the point is inside, while an even number or none places it outside. However, with a time complexity of $\mathcal{O}(n)$, ray tracing becomes inefficient for polygons with many vertices or when testing multiple points. It also struggles with edge cases, such as rays intersecting vertices or edges, requiring special handling. While it can handle non-convex polygons, its linear complexity is a bottleneck for large or irregular shapes. Additionally, it does not easily extend to 3D polyhedra, where feasibility testing requires more computationally intensive methods.

The computational time required for the R-function method is insignificant compared to other approaches, such as simplicial approximation and α-shape reconstruction.

Defining the feasible region analytically using R-functions provides several other advantages, such as facilitating operations like bounding box determination (boxing), and calculating the volume of the feasible region. One of the simplest methods to achieve these tasks is Monte Carlo sampling.

Steps for determining an Axis-Aligned Bounding Box:

1. Define the implicit function $R_\wedge$ that characterizes the feasible region.
2. Generate a set of points within a broader domain that fully encloses the feasible region.
3. Evaluate $R_\wedge$ at each sampled point and retain those where $R_\wedge \geq 0$.
4. Determine the minimum and maximum values of each parameter $x_k$ among the feasible points.

These values approximate the lower and upper bounds of a tightly fitting hyperrectangular box.

Another method for defining an axis-aligned bounding box involves solving a constrained nonlinear optimization problem:

$$\min_{x} \boldsymbol{\alpha}^\top \boldsymbol{x}, \quad \text{s.t.} \quad \|\boldsymbol{\alpha}\| = 1, \; R_\wedge(\boldsymbol{x}) = 0.$$

Minimizing $\boldsymbol{\alpha}^\top \boldsymbol{x}$ corresponds to projecting $\boldsymbol{x}$ linearly along direction $\boldsymbol{\alpha}$. By selecting different vectors $\boldsymbol{\alpha}$ and initial conditions for the optimization problem, one can determine the minimum value of $\boldsymbol{x}$ on the surface defined by $R_\wedge(\boldsymbol{x}) = 0$ along direction $\boldsymbol{\alpha}$.

## 7. Conclusions

In this work, we have demonstrated the effectiveness of V.L. Rvachev's R-functions in identifying and characterizing feasible regions in process design and operational flexibility analysis. R-functions provide a powerful algebraic framework for representing complex geometric regions and performing logical operations on constraint functions, allowing for an explicit analytical description of feasible domains.

Traditional feasibility analysis methods, such as optimization-based approaches and shape reconstruction techniques, often require extensive numerical computations, rely on convex approximations, or struggle to handle nonconvex and disjoint feasible regions. In contrast, the R-function methodology enables direct algebraic representation of feasibility conditions, eliminating the need for iterative numerical approximations and significantly improving computational efficiency.

We have applied the R-function methodology to multiple test cases, demonstrating its capability to define feasible regions analytically, even in cases involving nonlinear, nonconvex, and disjoint regions. Compared to conventional methods like simplicial approximation and α-shape reconstruction, the proposed approach offers several advantages, including the ability to handle implicit constraints, reduce computational complexity, and facilitate feasibility testing through a single analytical expression.

Compared to conventional methods like simplicial approximation and $\alpha$-shape reconstruction, the proposed approach eliminates the need for extensive sampling, provides an exact analytical representation of feasible regions, ensures precise boundary definition, and enables more efficient feasibility testing through direct evaluation of the created R-function instead of computationally expensive point-in-polygon tests.

When constraints are not explicitly defined by analytical functions but instead arise from complex models, simulations, or black-box functions, the methodology from Kucherenko et al. (2024) on analytical identification of Design Spaces can be applied. This involves approximating the constraints using a closed-form function, such as a multivariate polynomial model.

Moreover, our approach allows for additional operations such as bounding-box determination, volume computation, and feasibility classification using Monte Carlo methods or constrained optimization techniques. The use of R-functions is straightforward, computationally efficient, and scalable, as it can be applied with or without symbolic computation and is well-suited for high-dimensional problems. This flexibility makes R-functions particularly valuable in fields where defining feasible operational spaces is critical, such as pharmaceutical manufacturing, chemical process engineering, and energy systems.

## Acknowledgements

S.K and N.S. acknowledge the financial support of Eli Lilly and Company, the EPSRC Grant EP/T005556/1 and EPSRC grant EP/V042432/1.